\documentclass[11pt]{article}
\usepackage{epic,latexsym,amssymb}
\usepackage{amsfonts}
\usepackage{amscd}
\usepackage{amsmath}
\usepackage{graphicx}
\usepackage{color}
\usepackage{caption,
caption}
\usepackage{tikz}
\usepackage{amsthm}
\usepackage{bbm}
\usepackage{amsfonts}

\usepackage{caption}
\usepackage[margin=0.25in]{geometry}
\usepackage{pgfplots}
\pgfplotsset{width=10cm,compat=1.9}

\usepgfplotslibrary{external}
\usepackage{epic,latexsym,amssymb}
\usepackage{tikz}
\usepackage[bold,full]{complexity}
\usepackage{tikz,epsf}
\usepackage{hyperref}
\usepackage{graphicx}[width=\textwidth]
\usepackage{amsmath}
\usepackage{listings}

\usepackage{amsfonts}
\usepackage{amscd}
\usepackage{amsmath}
\usepackage{graphicx}
\usepackage{color}
\usepackage{caption,subcaption}
\usepackage{lineno}

\usepackage{graphicx}
\graphicspath{ {./images/} }

\usepackage[ruled,vlined]{algorithm2e}

\usepackage{lineno}

\newtheorem{thm}{Theorem}
\newtheorem{conj}[thm]{Conjecture}

\newcommand{\QEDmark}{\mbox{\textsc{qed}}}
\newcommand{\proofStarter}[1]{\textsc{#1} }

\def\vertex(#1){\put(#1){\circle*{2}}}
\def\vertexo(#1){\put(#1){\circle{2}}}
\def\vert(#1){\put(#1){\circle*{1.5}}}
\def\verto(#1){\put(#1){\circle{1.5}}}
\def\lab(#1)#2{\put(#1){\makebox(0,0)[c]{#2}}}
\definecolor{DarkGreen}{rgb}{0.2, 0.6, 0.3}

\definecolor{electricindigo}{rgb}{0.44, 0.0, 1.0}

\let\oldenumerate\enumerate
\renewcommand{\enumerate}{
  \oldenumerate
  \setlength{\itemsep}{0.5pt}
  \setlength{\parskip}{0pt}
  \setlength{\parsep}{0pt}
}

\begin{document}

\title{Advancements in Research Mathematics through AI: \\ A Framework for Conjecturing}
\author{$^{1, 2}$Randy Davila\
\\
\\
$^1$Research and Development \\
RelationalAI \\
Berkeley, CA 94704, USA\\
\small {\tt Email: randy.davila@relational.ai}\\
\\
$^2$Computational Applied Mathematics \\ \& Operations Research \\
Rice University \\
Houston, TX 77005, USA\\
\small {\tt Email: randy.r.davila@rice.edu}\\
\\
}

\date{}
\maketitle

\abstract{In the words of the esteemed mathematician Paul Erdös, the mathematician's task is to \emph{prove and conjecture}. These two processes form the bedrock of all mathematical endeavours, and in the recent years, the mathematical community has increasingly sought the assistance of computers to bolster these tasks. This paper is a testament to that pursuit; it presents a robust framework enabling a computer to automatically generate conjectures - particularly those conjectures that mathematicians might deem substantial and elegant. More specifically, we outline our framework and provide evidence in the mathematical literature demonstrating its use in generating publishable research and surprising mathematics. We suspect our simple description of computer-assisted mathematical conjecturing will catalyze further research into this area and encourage the development of more advanced techniques than the ones presented herein.}


\section{Introduction}\label{sec1}
In the grand edifice of mathematical knowledge, the questions we pose—our \emph{conjectures}—often serve as the keystone. These compelling, yet unresolved propositions delineate the boundaries of current understanding, and in doing so, chart the course of future discovery. They foster vibrant collaborations, prompt the development of novel mathematical disciplines and techniques, and continuously reinvigorate the field with fresh challenges.

Consider, for example, Fermat's Last Theorem. This assertion—that no three positive integers $a$, $b$, and $c$ can satisfy the equation $a^n + b^n = c^n$ for any integer value of $n > 2$—remained an open conjecture for nearly four centuries. Despite this, its very existence sparked the genesis of an entirely new mathematical discipline, algebraic number theory. Its resolution was finally accomplished by Andrew Wiles in 1994, a seminal achievement in the annals of mathematical history~\cite{Wiles}.

Similarly, the Four Color Theorem—the proposition that any planar map can be colored with just four colors such that no two adjacent regions share the same color—marked a watershed moment in the relationship between mathematics and computation. This conjecture, first posited in 1852, stood as a testament to the power of computer-assisted proof when it was resolved in 1976 with the aid of a computer program~\cite{4-color}.

The significance of conjectures is further underscored by the Clay Mathematics Institute's Millennium Prize Problems. These seven unsolved problems, each carrying a million-dollar prize for their resolution, have provided a north star for countless researchers worldwide. To date, only one has been solved—the Poincaré Conjecture by Grigori Perelman~\cite{Perelman1, Perelman2}. Yet the cumulative effort spent on these problems has precipitated remarkable advancements in diverse fields like weather forecasting, oil and gas extraction, network security, and quantum mechanics.

To attain such impact and popularity, a conjecture needs to be general, supported by evidence, simply stated, yet challenging to prove, or to disprove. Traditionally, professional mathematicians have formulated these conjectures, typically after months of building and analyzing examples of a particular mathematical object until a pattern or relationship is observed. This practice of conjecture-making has often been akin to a secret art form, with the source of the intuition sometimes carefully guarded. For instance, the legendary mathematician Gauss was known for his reticence in sharing the sources of his insights, often presenting elegant proofs that seemingly materialized `out of thin air,' leaving no trace of the exhaustive analysis and multitude of examples that underpinned his conclusions~\cite{Gauss}. In many ways, the proofs that we see published today mirror an iceberg—a refined text that can be comprehended in a few hours but represents the tip of an unseen mountain of months of trials, errors, and imperfect attempts, all suppressed for the sake of brevity.

Given that the process of conjecture-making involves the generation and analysis of relationships between various mathematical objects, it is logical to consider the potential role of computers in this domain. Indeed, artificial intelligence, machine learning, and deep learning have demonstrated their prowess in discerning hidden trends and relationships in vast data sets. As an example, a recent paper published in Nature~\cite{nature} demonstrated how these technologies can be harnessed to generate novel and impactful conjectures in the field of \emph{knot theory}. Yet, it is crucial to recognize that this is not a new phenomenon—the initial forays into computer-assisted conjecturing can be traced back to the late 1980s with Fajtlowicz's GRAFFITI program~\cite{Graffiti}. This pioneering work successfully generated significant (and valid) conjectures in graph theory, drawing the attention of esteemed mathematicians like Paul Erdös, Ronald Graham, and Fan Chung.

The domain of programming computers to formulate conjectures in mathematics is referred to as \emph{automated conjecturing}, an area pioneered by Fajtlowicz with the development of GRAFFITI. Many successors have followed in the footsteps of this groundbreaking work, including DeLaViña's GRAFFITI.pc~\cite{Graffitipc} and Larson's Conjecturing~\cite{Larson}, both of whom are notably graduate students of Fajtlowicz. There are also several other notable contributions in this domain, such as Lenat’s AM~\cite{Lenat_1, Lenat_2, Lenat_3}, Epstein’s GT~\cite{Epstein_1, Epstein_2}, Colton’s HR~\cite{Colton_1, Colton_2, Colton_3}, Hansen and Caporossi’s AGX~\cite{AGX_1, AGX_2, AGX_3}, Mélot’s Graphedron~\cite{graphedron_1}, and Davila's \emph{TxGraffiti}~\cite{TxGraffiti} and \emph{Conjecturing.jl}~\cite{Conjecturing.jl}.

These programs primarily adhere to a similar modus operandi: they calculate numerical properties on a given mathematical object, test for potential inequalities between these numerical properties, and then sieve the valid inequalities based on a heuristic. The heuristic is of paramount importance in this process. To comprehend why, it's crucial to acknowledge that a computer could churn out millions of relationships between parameters and properties for a given dataset. So, once a computer has generated a million potential relationships, the fundamental question that arises is, which subset of these relationships merits the distinction of a ``mathematical conjecture"?

In essence, the crux of automated conjecturing lies not merely in generating potential relationships, but in discerning which of these relationships have the potential to be mathematical conjectures. The answer to this question is complex, and it is determined by various factors, including the novelty of the conjecture, its complexity, and the potential impact it may have on the mathematical community. It is the heuristic of these programs that ultimately determines what subset of these relationships can be considered as ``interesting conjectures."

The \emph{Dalmatian} heuristic is one of the most prominent filtering methods, pioneered by Fajtlowicz in his GRAFFITI program. The foundation of Dalmatian is entrenched in the long-standing mathematician's philosophy of focusing only on the most significant conjectures, colloquially referred to as ``meaty conjectures." More specifically, suppose there is a set of computer-generated inequalities relating functions on a group of mathematical objects. If we have this set of inequalities, the Dalmatian heuristic would only present a new inequality as a conjecture if the new relationship achieved equality for at least one mathematical object in the database that did not attain equality for any previously stored relationships.

The strength of this heuristic lies in its ability to guarantee that no newly generated relationship would be a mere iteration of a known one. Hence, the conjectures generated by GRAFFITI and its successors were ``strong" in the sense that no other valid relationships ever surpassed them. This heuristic ensures that every newly proposed conjecture brings something unique and substantial to the table, preventing any potential redundancy in the generated conjectures. For a more in-depth explanation of this method, refer to~\cite{Larson}.

In this paper, we describe a novel automated conjecturing method implemented by the Python program TxGraffiti~\cite{TxGraffiti, TNP}, and the Julia programming package Conjecturing.jl~\cite{Conjecturing.jl}. We provide a detailed overview of the frameworks underlying these programs and illustrate their effectiveness with evidence from their application. 

\section{Framework Design}\label{sec:methods}

\subsection{Data Collection and Preparation}\label{subsec:data}
In the design of a computer program that generates mathematical conjectures, the first requirement is a database of mathematical objects. It is crucial to underscore the importance of data quality. An extensive database isn't necessary for the computer to identify non-trivial relationships among object properties. Rather, what is needed is a collection of unique instances of the objects in question, such as special counter-examples or interesting families of graphs from the literature. In our implementations, we utilized databases of several hundred objects, though we have also experimented with thousands. 

\subsection{Feature Generation}\label{subsec:feature}
\begin{figure}[h]
    \centering
    \includegraphics[width=11.5cm, height=3.5cm]{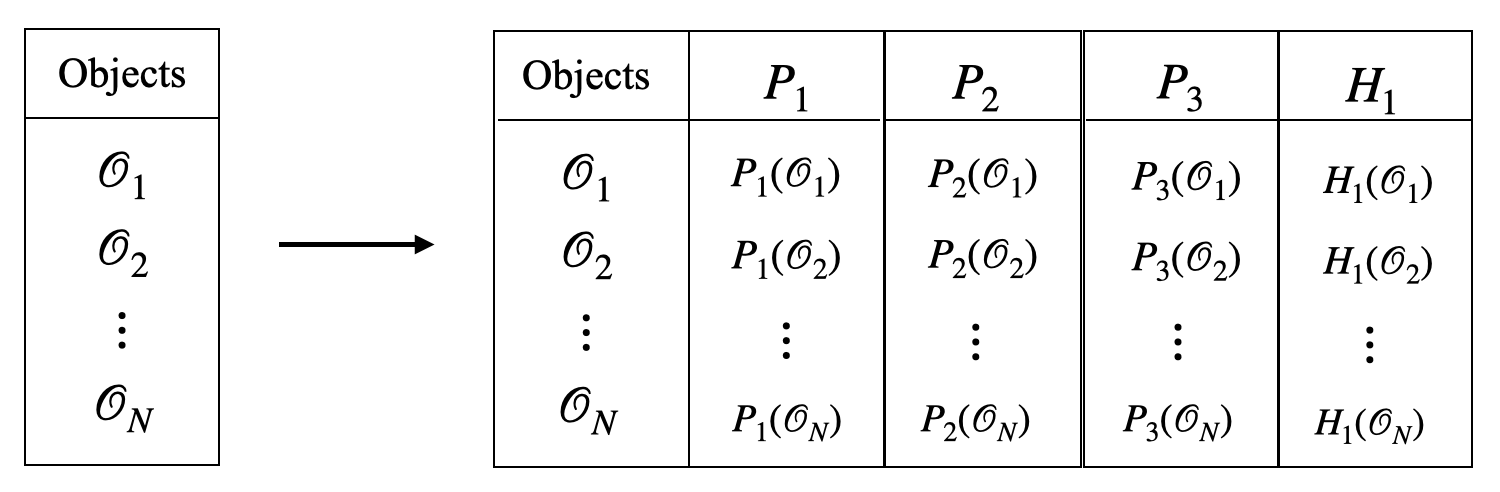}
    \caption{A mapping of a collection of $N$ mathematical objects to a table of numerical and Boolean properties.}
    \label{fig:1}
\end{figure}

After the database is stored on disk, the next step is to generate a table of various desired functions computed from the objects in this database. Our framework mandates that at least two of these functions return numerical values (for pairwise comparison), while others can return either numerical or Boolean values. See Figure~\ref{fig:1} for an illustration of this process; numerical properties are denoted by $P_i$, and Boolean properties are denoted by $H_i$. 

\subsection{Inequality Generation}\label{subsec:inequality}
We propose and implement a simplified version of the following steps for a computer program to generate inequalities relating properties of the objects under consideration, with an emphasis on simplicity and strength. 
\begin{enumerate}
    \item[1.] Select a target property $P_i$ -- a precomputed and numerically valued function on the objects in the database. 
    \item[2.] Choose an inequality direction (upper or lower) to bound the property $P_i$.
    \item[3.] For each precomputed numerical function $P_j$, with $j \neq i$, use a supervised machine learning technique or linear program to find a function $f$ such that $P_i(\mathcal{O}) \le f((P_j(\mathcal{O}))$ holds for each object $\mathcal{O}$ in the database, \emph{and} the number of instances where the inequality is an equality is maximized. 
    \item[4.] If $P_i(\mathcal{O}) \neq  f((P_j(\mathcal{O}))$ for all objects $\mathcal{O}$ in the database, disregard $f$ as a conjectured upper (or lower) bound on $P_i$. Otherwise, $f$ is called a \emph{sharp bounding function}; store $f(P_j)$ as a conjectured upper (or lower) bound on $P_i$ and record the set of objects $\mathcal{O}$ where  $P_i(\mathcal{O}) =  f(P_j(\mathcal{O}))$; the size of this set is the \emph{touch number} of the conjecture. 
\end{enumerate}

Both TxGraffiti and Conjecturing.jl follow these steps, with conjectured upper and lower bounds computed automatically through linear programming formulations. For example, consider producing one conjectured upper bound on $P_i(\mathcal{O})$ in terms of another numerical valued function $P_j(\mathcal{O})$). This is achieved by solving a linear programming problem. In the simplest case, TxGraffiti aims to minimize a some linear function $f(m, b)$ subject to a set of constraints.  
\begin{equation*}
\begin{aligned}
& \underset{m, b}{\text{minimize}}
& & f(m, b) \\
& \text{subject to}
& & P_i(\mathcal{O}) \leq mP_j(\mathcal{O}) + b, & \forall \mathcal{O} \in \text{Database},
\end{aligned}
\end{equation*}

\begin{figure}[h]
    \centering
    \includegraphics[width=10.5cm, height=4.25cm]{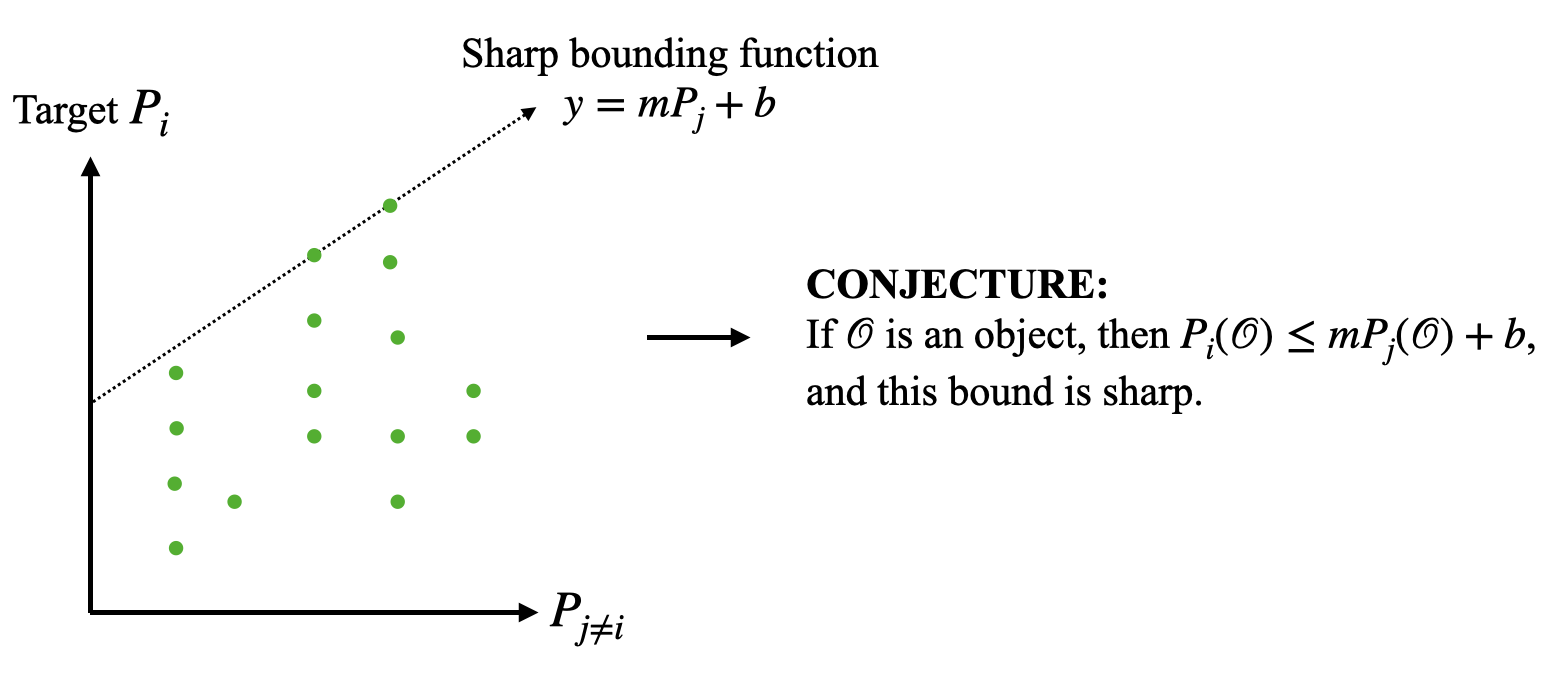}
    \caption{Finding a possible (linear) upper bound on the target property $P_i$ in terms of property $P_j$.}
    \label{fig:2}
\end{figure}

The goal is to find the line with the slope $m$ and y-intercept $b$ that satisfies all the inequalities and maximizes the number of times these inequalities hold with equality. In this way, TxGraffiti searches for the best linear upper bound on $P_i(\mathcal{O})$ in terms of $P_j(\mathcal{O})$ that holds for all objects $\mathcal{O}$ in the database; see Figure~\ref{fig:2} for a graphical illustration of the linear upper bound.

\begin{figure}[h]
    \centering
    \includegraphics[width=10.5cm, height=4.25cm]{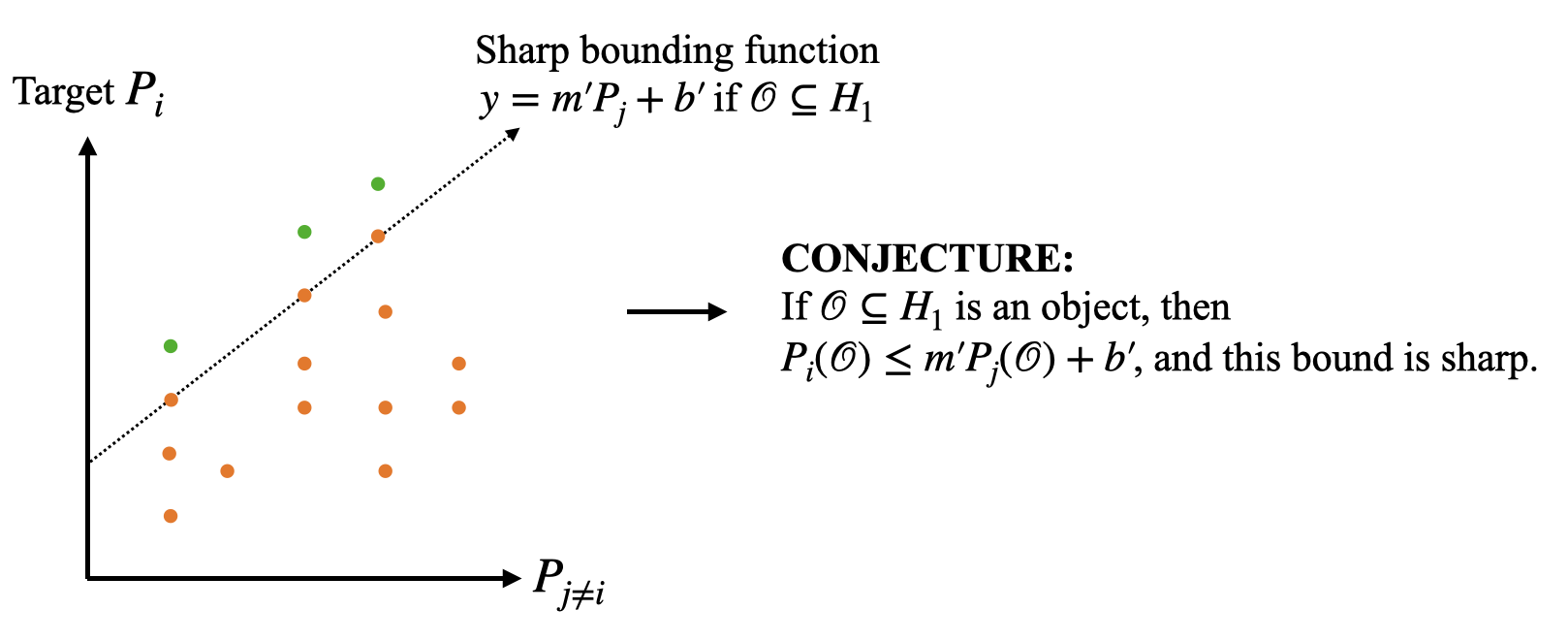}
    \caption{Finding a possible (linear) upper bound on the target property $P_i$ in terms of property $P_j$ for objects in $H_1$.}
    \label{fig:3}
\end{figure}

The above process is done on each numerical column of the feature data; that is, each pair of numerical functions is compared against each other and given a proposed inequality conjectured between them. Each conjecture generated by the above steps applies to all types of objects in the database. However, we can generate even more conjectures. By applying the same steps to a subset of objects in the database that satisfy a particular Boolean property (or combination of Boolean properties), we may obtain conjectures that are less general but potentially stronger; see Figure~\ref{fig:3} for an illustration of this process.

\subsection{Filtering and Sorting}\label{subsec:filter}

Once the initial steps are performed, we are left with a database of potential conjectures, along with detailed data for each conjecture. This data includes the set of objects that satisfy the conjecture's hypothesis, the graphs that attain equality, and the count of these graphs. At this point, our program implements a filtering heuristic to refine the list of conjectures.

Specifically, consider a target property $P_i$ and two conjectured inequalities $C_1$ and $C_2$, both providing upper (or lower) bounds on $P_i$. Suppose that the sharp bounding functions for $C_1$ and $C_2$ are identical and are both functions of a property $P_j$. If the set of objects that meet the hypothesis of $C_2$ is a strict subset of those that meet the hypothesis of $C_1$, then we remove $C_2$ from the database of stored conjectures.

This filtering process ensures that no stored conjecture is less general than another. In other words, we aim to prevent redundancy and to keep only the most general conjectures. It helps to streamline the database and ensures that each conjecture provides distinct and meaningful information.

Consider an illustrative example: suppose $C_1$ states that ``for all connected graphs, $\alpha(G) \leq 2\delta(G)$," and $C_2$ states that ``for all connected bipartite graphs, $\alpha(G) \leq 2\delta(G)$." Here, $\alpha(G)$ represents the independence number of $G$, and $\delta(G)$ represents the minimum degree $G$. Notice that the set of objects that satisfy the hypothesis of $C_2$ (connected bipartite graphs) is a strict subset of those that meet the hypothesis of $C_1$ (connected graphs). Therefore, $C_2$ would be removed during the filtering process, leaving only the more general $C_1$ in the stored database.

Finally, the conjectures stored are sorted according to their respective touch number. That is, conjectures of the program are presented to the user in non-increasing order concerning how many times the conjectured inequality holds with equality in the database of mathematical objects. After the conjectures have been sorted, the program shows the user a truncated list of the found inequalities; all posed as possible conjectures. 

\subsection{Code and Reproducibility}\label{subsec:code}
For readily available examples of this process, see the GitHub repositories~\cite{TxGraffiti, Conjecturing.jl}. The README file in each repository provides a step-by-step guide for running our code. 

\section{Results}\label{sec:results}
In this section, we highlight a series of substantial mathematical outcomes in the domain of \emph{graph theory}, spurred by conjectures derived from TxGraffiti and our framework. As an active and widely-utilized branch of mathematics that dates back to the era of the renowned mathematician Euler, graph theory centers around discrete mathematical structures that encapsulate relationships between vertices, manifested as connections or edges. For an in-depth exploration of graph theory and related terminology employed in this section, we recommend the seminal textbook by West~\cite{West}.

Our conjecturing framework, elaborated in the subsequent section, exhibits significant effectiveness when deployed on graphs, particularly those that are finite, simple, and undirected. To illustrate, when we tasked TxGraffiti to formulate conjectures on the extremal cardinalities of \emph{independent sets}, \emph{dominating sets}, and \emph{matching sets}, it generated a plethora of conjectures. This facilitated the development of several theorems connecting these graph invariants, resulting in outputs of a quality suitable for publication. Detailed presentations of these outcomes can be found in the publication by Caro et al.~\cite{CaDaHePe2022}. 

In a noteworthy instance, TxGraffiti conjectured that the independence number of an $r$-regular graph, where $r > 0$ (an $r$-regular graph is a graph wherein each vertex is linked to exactly $r$ other vertices by an edge), is at most the matching number. The original conjecture was framed for 3-regular graphs, a subset of $r$-regular graphs. A subsequent version of TxGraffiti extended the conjecture to encompass all $r$-regular graphs. This conjecture was later proven by Caro, Davila, and Pepper in~\cite{CaDaPe2020}, with Larson independently verifying the result.
\begin{conj}[TxGraffiti - Confirmed - Caro, Davila, and Pepper~\cite{CaDaPe2020}]\label{conj:1}
If $G$ is a connected $r$-regular graph with $r > 0$, independence number $\alpha(G)$, and matching number $\mu(G)$, then
\[
\alpha(G) \le \mu(G),
\]
and this bound is sharp. 
\end{conj}

The surprising nature of Conjecture~\ref{conj:1} is worth considering. Regular graphs, independence, and matchings are among the oldest and most deeply explored properties in graph theory. Interestingly, several mathematicians (based on correspondence with the author) initially perceived Conjecture~\ref{conj:1} as either obviously incorrect or known. However, upon thorough investigation, it became clear that Conjecture~\ref{conj:1} had not previously been proposed. Furthermore, it's noteworthy that the computation of the independence number is NP-hard, whereas the computation of the matching number is polynomial time. Therefore, the theorem resulting from Conjecture~\ref{conj:1} provides a tight upper bound on an NP-complete graph invariant in terms of a polynomial time graph invariant, rendering this result not only theoretically intriguing but also practically applicable.
 
Conjecture~\ref{conj:1} has spurred several other results relating matching sets and independent sets, which are also outlined in~\cite{CaDaPe2020}. TxGraffiti's conjecture about \emph{minimum maximal independent sets} versus \emph{minimum maximal matchings}, which remains open, is one such notable example.
\begin{conj}[TxGraffiti - Open]\label{conj:2}
If $G$ is a connected $r$-regular graph with $r > 0$, $i(G)$ denotes the cardinality of a minimum maximal independent set in $G$, and $\mu^*(G)$ denotes the cardinality of a minimum maximal matching in $G$, then
\[
i(G) \le \mu^*(G),
\]
and this bound is sharp. 
\end{conj}

The \emph{zero forcing number} is an intriguing invariant that is linked with the positioning and quantity of power monitoring units in electrical grids, and with the theoretical minimum rank problem of symmetric matrices~\cite{AIM-Workshop}. Unlike the independence number and the domination number, the zero forcing number acts like an infection, spreading from one vertex to another through a color change rule.

TxGraffiti has proposed significant conjectures on the zero forcing number of cubic graphs (graphs in which every vertex is connected to exactly three other vertices). These conjectures and their results are extensively covered in Davila's Dissertation~\cite{Davila2}. TxGraffiti made the following three conjectures. The first two were verified by Davila and Henning in~\cite{DaHe19b} and~\cite{DaHe21a}, respectively. The last conjecture is still open, although Davila and Henning were able to prove a limited version of it in~\cite{DaHe19c}; see Theorem~\ref{thm:1} below.
\begin{conj}[TxGraffiti - Confirmed - Davila and Henning~\cite{DaHe19b}]\label{conj:3}
If $G \neq K_4$ is a connected and cubic graph with zero forcing number $Z(G)$ and total domination number $\gamma_t(G)$, then
\[
Z(G)\le \frac{3}{2}\gamma_t(G),
\]
and this bound is sharp. 
\end{conj}

\begin{conj}[TxGraffiti - Confirmed - Davila and Henning~\cite{DaHe21a}]\label{conj:4}
If $G \neq K_4$ is a connected and cubic graph with zero forcing number $Z(G)$ and domination number $\gamma(G)$, then
\[
Z(G)\le 2\gamma(G),
\]
and this bound is sharp. 
\end{conj}

\begin{conj}[TxGraffiti - Open]\label{conj:main}
If $G \neq K_4$ is a connected graph with maximum degree $\Delta(G) \leq 3$, zero forcing number $Z(G)$ and independence number $\alpha(G)$, then
\[
Z(G)\le \alpha(G) + 1,
\]
and this bound is sharp. 
\end{conj}

\begin{thm}[Davila and Henning~\cite{DaHe19c}]\label{thm:1}
If $G$ is a connected, claw-free, and cubic graph on $n \geq 6$ vertices, and with zero forcing number $Z(G)$ and independence number $\alpha(G)$, then
\[
Z(G)\le \alpha(G) + 1,
\]
and this bound is sharp. 
\end{thm}

Interestingly, the independence number of a claw-free graph can be computed in polynomial time~\cite{Stable-NP}, while the complexity of computing the zero forcing number of a claw-free graph is unknown (it is generally NP-hard~\cite{NP-Complete}). Thus, Theorem~\ref{thm:1} provides valuable insights into the computational complexity of zero forcing in claw-free cubic graphs.

TxGraffiti's conjectures often give rise to more general statements. For instance, Conjecture~\ref{conj:5} below, which initially applied to claw-free graphs, was proven by Brimkov et al. in~\cite{TxGraffiti-2023}. The proof, which was algorithmic in nature, led to the discovery of the more general Theorem~\ref{thm:2} shown below.
\begin{conj}[TxGraffiti - Confirmed - Brimkov et al.~\cite{TxGraffiti-2023}]\label{conj:5}
If $G$ is a claw-free graph with forcing number $Z(G)$ and vertex cover number $\beta(G)$, then
\[
Z(G)\le \beta(G),
\]
and this bound is sharp. 
\end{conj}

\begin{thm}[Brimkov et al.~\cite{TxGraffiti-2023}]\label{thm:2}
If $G$ is a connected graph with maximum degree $\Delta \geq 3$, zero forcing number $Z(G)$, and vertex cover number $\beta(G)$, then 
$$
Z(G) \le (\Delta - 2)\beta(G) + 1,
$$
and this bound is sharp.
\end{thm}

\section{Discussion}\label{sec:discussion}
In this paper, we have elucidated our general approach to utilizing computer programming for generating novel and meaningful mathematical conjectures. We have demonstrated the validity and effectiveness of our framework with a number of original conjectures it has successfully generated. Particularly noteworthy is the discovery of a previously unknown correlation between two extensively studied graph theoretic parameters.

However, the implemented versions of our framework were deliberately simplified, primarily focusing on identifying linear relationships as sharp bounding functions. While this provides an important baseline for comparison, the true potential of the framework may not yet be fully realized.

In the realm of further research, more advanced machine learning techniques could be employed to extend our framework, potentially uncovering sharp bounding functions involving multiple variables. Our current methods were not designed to identify bounding functions of polynomial degree higher than one; exploring techniques to tackle these scenarios may open up new opportunities for more complex and nuanced conjectures.

Moreover, beyond the mathematical and algorithmic considerations, we must also grapple with philosophical and practical questions. How can we balance the desire for generality with the need for specificity in conjecture generation? How should we approach the challenge of confirming or refuting the conjectures that are produced? How can we best integrate these tools into the broader mathematical community?

While our framework has proven to be a valuable tool in stimulating mathematical discovery, we must continue to refine and expand it in response to these ongoing questions. The objective is not merely to produce a greater quantity of conjectures, but to enhance the quality and impact of the conjectures generated. As such, our work represents not a final destination, but an exciting point of departure for future explorations in the realm of automated conjecture generation. \\

\noindent\textbf{Acknowledgements} \\ 
We would like to take a moment to express our profound gratitude to the many mathematicians whose significant contributions have been instrumental to the development and advancement of TxGraffiti and Conjecturing.jl. Firstly, we extend our sincere thanks go to Siemion Fajtlowicz and Ermelinda DeLaViña, who initially ignited the spark that inspired the author to delve into the development of TxGraffiti and Conjecturing.jl. Next we extend our heartfelt gratitude to Yair Caro, Michael Henning, Houston Schuerger, Ryan Pepper, and Michael Young. Their persistent and meticulous efforts in identifying a multitude of counter-examples and generalizations significantly refined our conjecturing programs.

We are particularly grateful to David Amos for his invaluable assistance since the inception of TxGraffiti and Conjecturing.jl. His expertise in programming and his understanding of the mathematical principles at the heart of our software were indispensable in translating our ambitious ideas into reality. Our thanks also go to Boris Brimkov for his contribution to the exposition in this paper. 



\medskip
\begin{thebibliography}{99}


\bibitem{AGX_3} M. Aouchiche, G. Caporossi, P. Hansen, and M. Laffay, Autographix: A survey, \emph{Electron. Notes Discrete Math.} \textbf{22} (2005), 515--520.

\bibitem{AIM-Workshop} AIM Special Work Group, Zero forcing sets and the minimum rank of graphs, \textit{Linear Algebra Appl.}, \textbf{428 (7)} (2008), 1628--1648.

\bibitem{4-color} K. Appel, W. Haken, Every planar map is four colorable, \emph{Part I. Discharging
Illinois J. Math.}, \textbf{21} (1977), 429--490.

\bibitem{TxGraffiti-2023} B. Brimkov, R. Davila, H. Schuerger, and M. Young, Computer assisted discovery: Zero forcing vs vertex cover, available at \href{https://arxiv.org/pdf/2209.04552.pdf}{https://arxiv.org/pdf/2209.04552.pdf}, (2022). 

\bibitem{AGX_1} G. Caporossi and P. Hansen, Variable neighborhood search for extremal graphs: 1 The autographix system, \emph{Discrete Math.} \textbf{212} (1--2) (2000), 29--44.

\bibitem{AGX_2} G. Caporossi and P. Hansen, Variable neighborhood search for extremal graphs: 5 Three ways to automate finding conjectures, \emph{Discrete Math.} \textbf{276} (1--3) (2004), 81--94.

\bibitem{CaDaPe2020} Y. Caro, R. Davila, and R. Pepper, New results relating matching and independence, \textit{Discuss. Math. Graph Theory} \textbf{42} (2020), 921--935. 

\bibitem{CaDaHePe2022} Y. Caro, R. Davila, M.A. Henning, and R. Pepper, Conjectures of TxGraffiti: Independence, domination, and matchings, \emph{Australas. J. Comb.} \textbf{84} (2) (2022), 258--274. 

\bibitem{NP-Complete} C. Chekuri and N. Korula, A graph reduction step preserving element-connectivity and applications, \emph{Automata, Languages, and Programming,} (2009) 254--265. 

\bibitem{Colton_1} S. Colton, A. Bundy and T. Walsh, Automated concept formation in pure mathematics, \emph{Proc. of the 16th Int. Jt. Conf. on Artif. Intell.,} vol. 2, IJCAI'99, Morgan Kaufmann Publishers (1999), 786--791.

\bibitem{Colton_2} S. Colton, Refactorable numbers---a machine invention, \emph{J. Integer Seq.} \textbf{2} (1999), Article 99.1.2.

\bibitem{Colton_3} S. Colton, Automated Theory Formation in Pure Mathematics, Springer, Heidelberg (2002).


\bibitem{TxGraffiti} R. Davila, LinearTxGraffiti, available at \\ \href{https://github.com/RandyRDavila/Linear\_TxGraffiti}{https://github.com/RandyRDavila/Linear\_TxGraffiti} (2019).

\bibitem{TNP} R. Davila, Conjecturing.jl, available at \\ \href{https://github.com/RandyRDavila/TNP-That-New-Program}{https://github.com/RandyRDavila/TNP-That-New-Program} (2020).

\bibitem{Conjecturing.jl} R. Davila, Conjecturing.jl, available at \\ \href{https://github.com/RandyRDavila/Conjecturing.jl}{https://github.com/RandyRDavila/Conjecturing.jl} (2022).

\bibitem{Davila2} R. Davila, Total and Zero Forcing in Graphs, \emph{University of Johannesburg}, PhD Thesis (2019).

\bibitem{DaHe19b} R. Davila and M. A. Henning, Total forcing versus total domination in cubic graphs, \textit{Appl. Math. Comput.}, \textbf{354} (2019), 385--395.


\bibitem{DaHe19c} R. Davila and M.A. Henning, Zero forcing in claw-free cubic graphs, \emph{Bull. Malays. Math. Sci. Soc.}, \textbf{43}  (2020), 673--688. 

\bibitem{DaHe21a} R. Davila and M. A. Henning, Zero forcing versus domination in cubic graphs, \emph{J. Comb. Optim.}, \textbf{41} (2021), 553--577.

\bibitem{nature} A. Davies, P. Veličković, L. Buesing, S. Blackwell, D. Zheng, N. Tomašev, R. Tanburn, P. Battaglia, C. Blundell, A. Juhász, M. Lackenby, G. Williamson, D. Hassabis, and P. Kohli, Advancing mathematics by guiding human intuition with AI, \emph{Nature}, \textbf{600} (2021), 70--74.

\bibitem{Graffitipc} E. DeLaViña, Graffiti.pc: A variant of Graffiti, \emph{DIMACS Ser. Discret. Math. Theor. Comput. Sci.} \textbf{69} (2005), p. 71.

\bibitem{Graffiti} E. DeLaViña, Some history of the development of Graffiti, \emph{Graphs and Discovery, DIMACS Ser. Discret. Math. Theor. Comput. Sci.} \textbf{69}, Amer. Math. Soc., Providence, RI (2005), 81--118.


\bibitem{Epstein_1} S.L. Epstein, On the discovery of mathematical theorems, \emph{IJCAI} (1987), 194--197.

\bibitem{Epstein_2} S.L. Epstein, Learning and discovery: One system's search for mathematical knowledge, \emph{Comput. Intell.} \textbf{4} (1) (1988), 42--53.

\bibitem{Gauss} C. F. Gauss, Disquisitiones arithmeticae, (1801).

\bibitem{Larson} C. E. Larson and N. Van Cleemput. Automated conjecturing I: Fajtlowicz’s
Dalmatian heuristic revisited, \emph{Artif. Intell.} \textbf{231} (2016), 17--38. 

\bibitem{Lenat_1} D.B. Lenat, The ubiquity of discovery, \emph{Artif. Intell.} \textbf{9} (3) (1977), 257--285.

\bibitem{Lenat_2} D.B. Lenat, On automated scientific theory formation: A case study using the am program, \emph{Mach. Intell.} \textbf{9} (1979), 251--286.

\bibitem{Lenat_3} D.B. Lenat, The nature of heuristics, \emph{Artif. Intell.} \textbf{19} (2) (1982), 189--249.

\bibitem{claw-free} E. Flandrin, R. Faudree, and Z. Ryj{\'a}{\v{c}}ek, Claw-free graph - a survey, \emph{Discrete Math.} \textbf{214} (2016), 196--200. 

\bibitem{Stable-NP} D. Nakamura and A. Tamura, A Revision of Minty's Algorithm for Finding a Maximum Weight Stable Set of a Claw-Free Graph, \emph{J. Oper. Res. Soc. Japan} \textbf{44} (2001), 194--204.

\bibitem{graphedron_1} H. Mélot, Facet defining inequalities among graph invariants: The system graphedron, \emph{Discrete Appl. Math.} \textbf{156} (10) (2008), 1875--1891.

\bibitem{Perelman1} G. Perelman, The entropy formula for the Ricci flow and its geometric applications. Preprint at arXiv:math.DG/0211159.

\bibitem{Perelman2} G. Perelman, Ricci flow with surgery on three-manifolds. Preprint at arXiv:math.DG/ 0303109. 

\bibitem{AlphaGo} D. Silver, A. Huang, C. J. Maddison, A. Guez, L. Sifre, G. van den Driessche, J. Schrittwieser, I. Antonoglou, V. Panneershelvam, M. Lanctot, S. Dieleman, D. Grewe, J. Nham, N. Kalchbrenner, I. Sutskever, T. Lillicrap, M. Leach, K. Kavukcuoglu, T. Graepel, and D. Hassabis, Mastering the game of Go with deep neural networks and tree search, \emph{Nature} \textbf{529} (2016), 484--489.

\bibitem{protein} A. W. Senior, R. Evans, J. Jumper, J. Kirkpatrick, L. Sifre, T. Green, C. Qin, A. Žídek, A. W. R. Nelson, A. Bridgland, H. Penedones, S. Petersen, K. Simonyan, S. Crossan, P. Kohli, D. T. Jones, D. Silver, K. Kavukcuoglu, and D. Hassabis, Improved protein structure prediction using potentials from deep learning, \emph{Nature} \textbf{577} (2020), 706--710. 

\bibitem{Turing} A. Turing, Intelligent machinery. \emph{The Essential Turing}, (2004), 395--432.

\bibitem{West} D. B. West, Introduction to Graph Theory 2nd Edition. Prentice-Hall (20010. ISBN: 0-13-014400-2 (print)

\bibitem{Wiles} A. Wiles, Modular Elliptic Curves and Fermat’s Last Theorem, \emph{Annals of Mathematics}, \textbf{141} (1995), 443--551.

\end{thebibliography}
\end{document}